# STEP-UP SIMULTANEOUS TESTS FOR IDENTIFYING ACTIVE EFFECTS IN ORTHOGONAL SATURATED DESIGNS

By Samuel S. Wu and Weizhen Wang[1]

*University of Florida and Wright State University*

A sequence of null hypotheses regarding the number of negligible effects (zero effects) in orthogonal saturated designs is formulated. Two step-up simultaneous testing procedures are proposed to identify active effects (nonzero effects) under the commonly used assumption of effect sparsity. It is shown that each procedure controls the experimentwise error rate at a given $\alpha$ level in the strong sense.

**1. Introduction.** Assume a linear model

$$(1) \qquad Y_i = \mu + \beta_1 x_{i1} + \cdots + \beta_k x_{ik} + \varepsilon_i, \qquad \text{for } i = 1, \ldots, M,$$

where $\varepsilon_i \sim$ i.i.d. $N(0, \sigma^2)$. The unknown parameters $\beta_i$ are of interest and $\mu$ and $\sigma$ are two unknown nuisance parameters. The design is called *orthogonal* if the least squares estimators $\hat{\beta}_i$ ($1 \leq i \leq k$) are uncorrelated (equivalent to independent), which occurs, for example, in two-level fractional factorial designs. The design is said to be *saturated* if there are just enough observations to estimate the model parameters $\beta_i$ and $\mu$ (i.e., $M = k + 1$), leaving no degrees of freedom to estimate the error variance $\sigma^2$. In order to make inferences on $\beta_i$, one must typically use the assumption of effect sparsity, that is, that most of the $\beta_i$'s are equal to zero. Then we can use the corresponding $\hat{\beta}_i$'s to estimate $\sigma^2$. However, we do not know how many and which of the $\beta_i$'s are zero. An initial guess would be at least $\nu$ of the $\beta_i$'s equal zero, say one-half of the effects. Therefore, the smallest $\nu$ of the $\hat{\beta}_i^2$'s should be used to estimate $\sigma^2$. Any other $\hat{\beta}_j$ whose square is substantially larger is likely to have a nonzero mean and corresponds to an active effect.

For a fixed sequence of $\underline{\beta} = (\beta_1, \ldots, \beta_k)$, let

$$(2) \qquad N = \text{the number of } \beta_i\text{'s which equal zero}.$$

Received May 2004; revised October 2005.
[1]Supported in part by NSF Grant DMS-03-08861.

Thus, the number of nonzero $\beta_i$'s is equal to $k - N$ and the entire parameter space without nuisance parameters is $H = \{\underline{\beta} = (\beta_1, \ldots, \beta_k) : N \geq \nu\}$. For each integer $m \in [\nu + 1, k]$, consider the testing problem

(3) $\quad H_{0,m} : N \geq m \quad \text{vs.} \quad H_{A,m} : N \leq m - 1 \quad$ (i.e., $k - N \geq k - m + 1$)

and define a parameter configuration in each $H_{0,m}$,

(4) $$\underline{\beta}_m =: (0, \ldots, 0, +\infty, \ldots, +\infty),$$

where the first $m$ components are zero. Let

(5) $$\mathcal{B} = \{H_{0,m} : \nu + 1 \leq m \leq k\},$$

which contains all null hypotheses of interest in this paper. Because $H_{0,i}$ is a subset of $H_{0,j}$ for any $i > j$, if $H_{0,j}$ is incorrect, then so is $H_{0,i}$. This implies that a testing process should be terminated as soon as a rejection occurs for some null hypothesis. Starting from $m = \nu + 1$, we test these hypotheses one at a time as $m$ goes up to $k$. If $H_{0,\nu+1}$ is rejected, we then conclude that there are $k - \nu$ active effects (i.e., $H \cap H_{A,\nu+1}$) and no longer test any other hypotheses; otherwise, we test the next hypothesis $H_{0,\nu+2}$. In general, if $H_{0,m_0}$ is the first hypothesis being rejected for some $m_0 \leq k$, we stop and conclude that there are $k - m_0 + 1$ nonzero effects (i.e., $H_{0,m_0-1} \cap H_{A,m_0}$); otherwise, all hypotheses in $\mathcal{B}$ are accepted and we conclude that there is no active effect. Clearly, this is a step-up testing procedure.

Many inference procedures have been proposed to identify active effects. The data analysis of orthogonal saturated designs was initially considered by Birnbaum [1] and Daniel [2]. The half-normal plot introduced by Daniel [2] is still being used in the preliminary analysis. Lenth [8] proposed the first adaptive method to let the data determine which and how many of the $\hat{\beta}_i$'s should be used to estimate $\sigma^2$. Whether Lenth's interval is of level $1 - \alpha$ still remains a question. Besides using the data adaptively, another fundamentally desirable property is the ability to control the error rate in the strong sense (i.e., under all parameter configurations), which is espoused by Hochberg and Tamhane [5]. For orthogonal saturated designs, the first adaptive confidence interval known to provide strong control of error rates and more general results can be found in [15] and [16], respectively. Hamada and Balakrishnan [4] provided a thorough review of the analysis methods available for saturated designs.

Since we do not know which and how many effects are active, it is reasonable to search for active effects using simultaneous tests. Here are two possibilities:

(a) Starting from the largest $\hat{\beta}_i^2$, test whether the corresponding effect is active. If it is not active, then conclude no active effect and stop; otherwise, test for the second largest $\hat{\beta}_i^2$ (go down), and so on, until a zero effect is found—the step-down tests.

(b) Starting from the $(\nu+1)$th smallest $\hat{\beta}_i^2$, test whether the corresponding effect is active. If it is active, then conclude $k-\nu$ active effects and stop; otherwise, test for the $(\nu+2)$th smallest $\hat{\beta}_i^2$ (go up), and so on, until an active effect is found—the step-up tests.

Voss [12] proposed nonadaptive step-down tests for a set of hypotheses different from $\mathcal{B}$ and controlled the experimentwise error rates at the given level $\alpha$ in the strong sense. Recently, Voss and Wang [14] derived adaptive step-down tests with the experimentwise error rates controlled in the same setting as Voss [12]. As pointed out by several researchers [3, 7, 11], the step-up tests typically have a greater power to detect active effects than the step-down tests. For example, Venter and Steel [11] proposed one step-up and one step-down procedure for the orthogonal saturated designs with cutoff points determined at $\underline{\boldsymbol{\beta}}_m$. However, they were unable to prove that the experimentwise error rates of their procedures are controlled at a given level $\alpha$ in the strong sense.

The rest of this article is organized as follows. In Section 2, we provide the motivation for two desirable tests for $H_{0,m}$ and establish a probability inequality for noncentral $\chi^2$ distributions which asserts that the maximal type I error of any test in a general class is always achieved at $\underline{\boldsymbol{\beta}}_m$. Based on this inequality, two level-$\alpha$ tests are proposed in Section 3 for testing each single hypothesis $H_{0,m}$. Two sequential step-up procedures for testing all hypotheses in $\mathcal{B}$, which control the experimentwise error rates at $\alpha$ in the strong sense, are derived in Section 4. Section 5 presents a simulation study and Section 6 concludes with some discussion.

**2. Motivation and a general probability inequality.** In this section, we provide motivation for the new procedures testing $H_{0,m}$ and a general class [see $\mathcal{A}_m$ in (11)] of the rejection regions of level-$\alpha$. The maximal type I error of each rejection region in this class is always achieved at $\underline{\boldsymbol{\beta}}_m$, as stated in Theorem 1.

Assume the factorial effect estimators $\hat{\beta}_i$ are independently distributed as $N(\beta_i, a_i^2 \sigma^2)$ for known constants $a_i$. We may assume that each $a_i = 1$ without loss of generality. Let $X_1, \ldots, X_k$ be the order statistics of the $\hat{\beta}_i^2$ for $i \leq k$. Intuitively, it is more likely that the small order statistics $X_i$ will correspond to estimators $\hat{\beta}_j$ with $\beta_j = 0$. If we believe a priori that at least $\nu$ of the $\beta_i$'s are zero, then those $\beta_i$'s corresponding to $X_1, \ldots, X_\nu$ are likely to be the negligible ones. To test $H_{0,m}$, one needs to compare $X_m$ with $X_1, \ldots, X_\nu$. For any integer $n \in [\nu, m]$, let $S_n = \sum_{i=1}^n X_i$ and $\bar{X}_n = S_n/n$ and define a test statistic for $H_{0,m}$ as follows:

$$(6) \qquad W_{n,m} = \frac{nX_m}{\sum_{i=1}^n X_i} = \frac{nX_m}{S_n} = \frac{X_m}{\bar{X}_n}.$$

Intuitively, a large value of $W_{n,m}$ is in favor of $H_{A,m}$. Therefore, the rejection region should be $W_{n,m} > d_{n,m}$ for a constant $d_{n,m}$ which satisfies

$$\sup_{\underline{\beta} \in H_{0,m}} P_{\underline{\beta}}(W_{n,m} > d_{n,m}) = \alpha. \tag{7}$$

In this section, we will show that $W_{n,m}$ is stochastically largest at $\underline{\beta}_m$ for $\underline{\beta} \in H_{0,m}$.

DEFINITION 1.  A random variable $X$ is said to be stochastically smaller than $Y$, denoted by $X \prec Y$, if $P(X \leq d) \geq P(Y \leq d)$ for all $d$.

Let $Y_{1,m}, \ldots, Y_{m,m}$ be the order statistics of the $\hat{\beta}_i^2$'s for $i \leq m$ when $\beta_1 = \cdots = \beta_m = 0$ and let

$$Z_{n,m} = \frac{nY_{m,m}}{\sum_{i=1}^{n} Y_{i,m}}. \tag{8}$$

Langsrud and Næs [7] also studied $Z_{n,m}$, called $\Psi_{m,0,n,m}$-distribution in their notation. Clearly, the distribution of $Z_{n,m}$ does not depend on any parameter and can be sampled based on the order statistics of $m$ independent $\chi_1^2$ random variables. It is easy to see that $W_{n,m} = Z_{n,m}$ at $\underline{\beta}_m$ and we want to show that $W_{n,m} \prec Z_{n,m}$ on $H_{0,m}$. If this is true, then (7) reduces to

$$P(Z_{n,m} > d_{n,m}) = \alpha \tag{9}$$

and $d_{n,m}$ is the $100(1-\alpha)$ percentile of $Z_{n,m}$.

DEFINITION 2.  A function $h: \mathbb{R}^d \to \mathbb{R}$ will be called nondecreasing (to the coordinatewise ordering) if $x_i \leq y_i$, $i = 1, \ldots, d$, implies that $h(\boldsymbol{x}) \leq h(\boldsymbol{y})$.

Now we prove a general theorem that includes (9) as a special case.

THEOREM 1.  *Suppose that $X_1, \ldots, X_k$ are the order statistics of independent random variables with noncentral chi-squared distributions $\chi_1^2(\beta_i^2)$, $1 \leq i \leq k$. Then*

$$\sup_{\underline{\beta} \in H_{0,m}} P_{\underline{\beta}}(R) = P_{\underline{\beta}_m}(R) \tag{10}$$

*for any rejection region $R$ of $H_{0,m}$ that belongs to the class*

$$\mathcal{A}_m = \{T_L(X_1, \ldots, X_{s-1}) < T_R(X_s), \text{ for some integer } 1 < s \leq m\}, \tag{11}$$

*where $T_L$ and $T_R$ satisfy the following two properties:*

(i) *(monotone) $T_L(x_1, \ldots, x_{s-1})$ and $T_R(x_s)$ are nondecreasing functions;*

(ii) (*invariant*) $T_L/T_R$ is invariant to scale transformation, that is, for any $a > 0$,

$$\frac{T_L(ax_1,\ldots,ax_{s-1})}{T_R(ax_s)} = \frac{T_L(x_1,\ldots,x_{s-1})}{T_R(x_s)}. \tag{12}$$

COROLLARY 1. *For any $0 < \alpha < 1$ and $d_{n,m}$ given in (9), the rejection region*

$$R_{n,m} = \{W_{n,m} > d_{n,m}\} \tag{13}$$

*defines a level-$\alpha$ test for (3).*

PROOF. Let $s = m$, $T_L(x_1,\ldots,x_{m-1}) = \sum_{i=1}^n x_i/n$ and $T_R(x_m) = x_m/d_{n,m}$. The claim follows from Theorem 1. □

To prove Theorem 1, we need the facts below (Corollary 2 and Corollary 3).

LEMMA 1. *For a nonnegative random variable $X$ and a positive number $y$, let $X_y = X/y$, given that $X \leq y$. Let $U \sim \chi_1^2$, a chi-squared distribution with one degree of freedom, and $V \sim \chi_1^2(\theta^2)$, a noncentral chi-squared distribution with one degree of freedom and noncentrality parameter $\theta^2$. Then we have* (i) $U_y \prec V_y$ *and* (ii) $U_{y_2} \prec U_{y_1}$ *for any $y_2 > y_1 > 0$.*

Part (ii) of Lemma 1 is identical to Lemma 2 of [7]. The proofs for both stochastic orderings follow from the monotone likelihood ratio function.

LEMMA 2. *Let $U_1,\ldots,U_{s-1}$ be any independent random variables. Let the same be true for $V_1,\ldots,V_{s-1}$. If $T(x_1,\ldots,x_{s-1})$ is a nondecreasing function and $U_i \prec V_i$ for $i \leq s-1$, then*

$$T(U_1,\ldots,U_{s-1}) \prec T(V_1,\ldots,V_{s-1}).$$

This is called by some researchers (see, e.g., [13]) a "stochastic ordering lemma."

COROLLARY 2. *Suppose that $U_i \sim \chi_1^2$, $1 \leq i \leq s-1$, and $V_j \sim \chi_1^2(\theta_j^2)$, $1 \leq j \leq s-1$, are independent random variables. For any $y > 0$, let $U_{(i),y}$ be the order statistics of $U_{i,y}$ ($= U_i/y$, given that $U_i \leq y$, as defined in Lemma 1) and let $V_{(j),y}$ be the order statistics of $V_{j,y}$. Then for any nondecreasing function $T_L$,*

$$T_L(U_{(1),y},\ldots,U_{(s-1),y}) \prec T_L(V_{(1),y},\ldots,V_{(s-1),y}). \tag{14}$$

*Therefore, for any $t > 0$,*

$$P(T_L(U_{(1),y},\ldots,U_{(s-1),y}) \leq t) \geq P(T_L(V_{(1),y},\ldots,V_{(s-1),y}) \leq t).$$

Proof. Let
$$T_L^*(U_{1,y},\ldots,U_{s-1,y}) = T_L(U_{(1),y},\ldots,U_{(s-1),y}).$$

Since each $U_{(i),y}$ is a nondecreasing function in each $U_{j,y}$ and $T_L$ is nondecreasing in each of its arguments, $T_L^*$ is nondecreasing in each of its arguments. Therefore, if one combines part (i) of Lemma 1 with Lemma 2, it can be concluded that

$$T_L(U_{(1),y},\ldots,U_{(s-1),y}) = T_L^*(U_{1,y},\ldots,U_{s-1,y}) \prec T_L^*(V_{1,y},\ldots,V_{s-1,y})$$
$$= T_L(V_{(1),y},\ldots,V_{(s-1),y}). \qquad \square$$

COROLLARY 3. *Suppose that $T_L$ and $T_R$ satisfy the monotone and invariant conditions specified in the definition of $\mathcal{A}_m$ in* (11) *and one defines*

$$G^*(y) \equiv P(T_L(U_{(1),y},\ldots,U_{(s-1),y}) \le T_R(1))$$
$$= \int\cdots\int_{\{0<u_1<\cdots<u_{s-1}\le y, T_L(u_1,\cdots,u_{s-1})\le T_R(y)\}} (s-1)! \frac{\prod_{i=1}^{s-1} f(u_i)\,du_i}{F(y)^{s-1}}.$$

*Then under the condition of Corollary* 2, $T_L(U_{(1),y},\ldots,U_{(s-1),y})$ *is stochastically nonincreasing in $y$, that is, for $y_2 > y_1 > 0$,*

$$T_L(U_{(1),y_2},\ldots,U_{(s-1),y_2}) \prec T_L(U_{(1),y_1},\ldots,U_{(s-1),y_1}).$$

*Therefore, $G^*(y)$ is a nondecreasing function.*

PROOF. Combine part (ii) of Lemma 1 with Lemma 2 and define $T_L^*$ as in Corollary 2. $\square$

PROOF OF THEOREM 1. Consider two samples,
$$\{\hat{\beta}_1^2,\ldots,\hat{\beta}_m^2, +\infty,\ldots,+\infty \ (k-m \text{ of them})\} \quad \text{and} \quad \{\hat{\beta}_1^2,\ldots,\hat{\beta}_k^2\}.$$

Clearly, the $s$th order statistic of the first sample is stochastically larger than that of the second sample, that is, $X_s \prec Y_{s,m}$ for any given integer $s$ satisfying $1 < s \le m$.

Second, we have
$$P(T_L(Y_{1,m},\ldots,Y_{s-1,m}) < T_R(Y_{s,m}))$$
$$= E[P(T_L(Y_{1,m},\ldots,Y_{s-1,m}) < T_R(y)|Y_{s,m}=y)]$$
$$= E[P(T_L(U_{(1),y},\ldots,U_{(s-1),y}) < T_R(1)|Y_{s,m}=y)]$$
$$= E[G^*(Y_{s,m})].$$

Next, for any partition $\omega = (j_1, j_2, \ldots, j_{s-1})(j_s)(j_{s+1}, \ldots, j_k)$ of the integers 1 to $k$, we denote by $E_\omega$ the event $\{\hat{\beta}_{j_i}^2 < \hat{\beta}_{j_s}^2, \forall i < s;\ \hat{\beta}_{j_s}^2 < \hat{\beta}_{j_l}^2, \forall l > s\}$. We note that for any $\underline{\beta} \in H_{0,m}$,

$P_{\underline{\beta}}(R)$

$= E_{\underline{\beta}}[P_{\underline{\beta}}(T_L(X_1, \ldots, X_{s-1}) < T_R(y) | X_s = y)]$

$= E_{\underline{\beta}}[P_{\underline{\beta}}(T_L(X_{1,y}, \ldots, X_{s-1,y}) < T_R(1) | X_s = y)]$

$= E_{\underline{\beta}}\left[\sum_\omega P_{\underline{\beta}}(\{T_L(X_{1,y}, \ldots, X_{s-1,y}) < T_R(1)\} \cap E_\omega | X_s = y)\right]$

$= E_{\underline{\beta}}\left[\sum_\omega P_{\underline{\beta}}(\{T_L(X_{1,y}, \ldots, X_{s-1,y}) < T_R(1)\} | E_\omega, X_s = y) P_{\underline{\beta}}(E_\omega | X_s = y)\right]$

$\leq E_{\underline{\beta}}\left[\sum_\omega P_{\underline{\beta}_m}(\{T_L(X_{1,y}, \ldots, X_{s-1,y}) < T_R(1)\} | X_s = y) P_{\underline{\beta}}(E_\omega | X_s = y)\right]$

$= E_{\underline{\beta}}[P_{\underline{\beta}_m}(T_L(X_{1,y}, \ldots, X_{s-1,y}) < T_R(1) | X_s = y)]$

$= E_{\underline{\beta}}[P(T_L(U_{(1),y}, \ldots, U_{(s-1),y}) < T_R(1) | X_s = y)]$

$= E_{\underline{\beta}}[G^*(X_s)],$

where the above inequality is due to Corollary 2.

Finally, since $G^*(y)$ is a nondecreasing function due to Corollary 3, the inequality $X_s \prec Y_{s,m}$ implies that $G^*(X_s) \prec G^*(Y_{s,m})$. Therefore, it can be concluded that

$$E_{\underline{\beta}}[G^*(X_s)] \leq E[G^*(Y_{s,m})]. \qquad \square$$

**3. Two tests for $H_{0,m}$.** In this section, we present two tests for $H_{0,m}$. While both are of level $\alpha$, they are to be used under different circumstances. Let $d_{\nu,m}$ and $d_{m-1,m}$ denote the numbers determined by (9) when $n = \nu$ and $n = m - 1$, respectively.

THEOREM 2. *For any $0 < \alpha < 1$, the rejection regions*

(15) $\quad R_{\nu,m} = \{W_{\nu,m} > d_{\nu,m}\} \quad and \quad R_{m-1,m} = \{W_{m-1,m} > d_{m-1,m}\}$

*both define level-$\alpha$ tests for (3).*

This theorem is a special case of Corollary 1 when $n = \nu$ and $n = m - 1$. The test based on $R_{\nu,m}$ estimates the error variance by $\bar{X}_\nu = \sum_{i=1}^\nu X_i / \nu$, irrespective of the value of $m$, while the test based on $R_{m-1,m}$ uses the adaptive estimator $\bar{X}_{m-1} = \sum_{i=1}^{m-1} X_i / (m-1)$. These are referred to by Venter and Steel [11] as *fixed* and *sequential scaling*, respectively. Region $R_{\nu,m}$ should be applied only when $H_{0,m}$ is of interest and no information is available on whether $H_{0,n}$ is true for any $n < m$. In such a case, we are certain that at least $\nu$ of the $\hat{\beta}_i^2$'s have a zero mean. It is reasonable to compare $X_m$ with the average of $X_1, \ldots, X_\nu$, the smallest $\nu$ of the $\hat{\beta}_i^2$'s, through their average. If $W_{\nu,m}$ is large, then one concludes that all $\hat{\beta}_i^2$'s corresponding to $X_m, \ldots, X_k$ are from populations with nonzero means. On the other hand, if one tests $H_{0,n}$ for $n < m$ sequentially up to $H_{0,m}$ and $H_{0,m-1}$ is accepted, then at least $m-1$ of the $\beta_i$'s are zero. In this case, one should compare $X_m$ with $X_1, \ldots, X_{m-1}$ and a large value of $W_{m-1,m}$ would lead to a rejection of $H_{0,m}$.

**4. Two step-up testing procedures.** When effect sparsity is assumed, we do not know which and how many of the $\beta_i$'s are zero. It is of more interest to conduct tests simultaneously to identify the nonzero effects. The tests developed in the previous section, in fact, can detect whether there is a jump at $X_m$ among $X_1, \ldots, X_k$. However, these tests cannot tell whether the jump, if it exists, is the first one, which is what interests us. Therefore, as mentioned earlier, since $H_{0,m}$ decreases as $m$ increases, one needs to conduct tests sequentially. Like all testing problems, there are two major concerns: to control the experimentwise error rate at a given level $\alpha$, that is,

$$(16) \quad \sup_{\underline{\beta} \in \bigcup_{m=\nu+1}^{k} H_{0,m}} P_{\underline{\beta}}(\text{assert not } H_{0,n}, \text{ which contains } \underline{\beta}, \text{ for some } n \in [\nu+1, k]) \leq \alpha,$$

and to obtain more powerful tests, which means larger rejection regions.

The first requirement (16) can be ensured by using the closed test procedure proposed by Marcus, Peritz and Gabriel [9]. For details, see [6], page 137. A naïve solution is to assert not $H_{0,m}$ (i.e., to assert $H_{A,m}$) if one rejects $H_{0,i}$ at level $\alpha$ for all $i \geq m$. For example, suppose $R_{\nu,i}$ is used to test for each $H_{0,i}$. Then assert not $H_{0,m}$ iff $R_m = \bigcap_{i=m}^{k} R_{\nu,i}$ is true. This, by the closed test procedure, controls the experimentwise error rate at $\alpha$. Note that $R_m$ decreases as $m$ increases, which contradicts the fact that $H_{0,m}$ is decreasing (we need $R_m$ to increase). Therefore, simply applying the closed test procedure on the tests derived in the previous section only results in less powerful tests for the simultaneous hypotheses. We require that the rejection region for $H_{0,m}$ (a) increases as $m$ gets larger and (b) is of level-$\alpha$. In this section, two testing procedures are discussed with their rejection regions denoted by $\{R^*_{\nu,m}\}_{m=\nu+1}^{k}$ and $\{R^*_{m-1,m}\}_{m=\nu+1}^{k}$ corresponding to $R_{\nu,m}$ and $R_{m-1,m}$, respectively.

4.1. *The construction of $\{R^*_{\nu,m}\}^k_{m=\nu+1}$: step-up tests with fixed scaling (SUF).* The general form of $R^*_{\nu,m}$, for $\nu+1 \leq m \leq k$, is

$$
(17) \qquad R^*_{\nu,m} = \bigcup_{i=\nu+1}^{m} \{W_{\nu,i} > d^*_{\nu,i}\} = \{S_\nu < \max\{\nu X_i/d^*_{\nu,i}\}^m_{i=\nu+1}\},
$$

where the sequence of constants $\{d^*_{\nu,m}\}^k_{m=\nu+1}$ is determined iteratively below. More precisely, $d^*_{\nu,m}$ depends on $d^*_{\nu,i}$ for $i < m$ and causes $R^*_{\nu,m}$ to have level-$\alpha$. It is clear that $R^*_{\nu,m}$ is nondecreasing when $m$ gets larger and is strictly increasing if all $d^*_{\nu,m}$ are finite. We start with the following lemma.

LEMMA 3. *For a sequence of random variables $\{\Delta_i\}^s_{i=0}$, where $s$ is a given positive integer,*

$$
(18) \qquad P(\Delta_0 < \max\{\Delta_i\}^s_{i=1}) \leq \sum_{i=1}^{s} P(\max\{\Delta_j\}^{i-1}_{j=0} < \Delta_i).
$$

PROOF. We prove (18) by induction. When $s = 1$, (18) is true. Suppose (18) is true for any $s = n$. Then for $s = n+1$,

$$
P(\Delta_0 < \max\{\Delta_i\}^{n+1}_{i=1}) \leq P(\Delta_0 < \max\{\Delta_i\}^n_{i=1}) + P(\max\{\Delta_j\}^n_{j=0} \leq \Delta_0 < \Delta_{n+1})
$$
$$
\leq \sum_{i=1}^{n+1} P(\max\{\Delta_j\}^{i-1}_{j=0} < \Delta_i). \qquad \square
$$

We now determine the sequence of constants $\{d^*_{\nu,m}\}^k_{m=\nu+1}$ starting from $m = \nu + 1$. For testing $H_{0,\nu+1}$, let $R^*_{\nu,\nu+1} = R_{\nu,\nu+1}$. It is a level-$\alpha$ test by Theorem 2.

For any $\nu + 2 \leq m < k$, let $d^*_{\nu,\nu} = \infty$ and suppose that $\{d^*_{\nu,i}\}^{m-1}_{i=\nu+1}$ are available. Then $d^*_{\nu,m}$ is determined by solving

$$
\sum_{i=\nu+1}^{m} P_{\underline{\beta}_m}(A_i) = \alpha,
$$
(19)
$$
\text{where } A_i = \{\max\{S_\nu, \{\nu X_j/d^*_{\nu,j}\}^{i-1}_{j=\nu}\} < \nu X_i/d^*_{\nu,i}\}.
$$

Note that $R^*_{\nu,m}$ is of level-$\alpha$ because for any $\underline{\beta} \in H_{0,m}$,

$$
(20) \qquad P_{\underline{\beta}}(R^*_{\nu,m}) \leq \sum_{i=\nu+1}^{m} P_{\underline{\beta}}(A_i) \leq \sum_{i=\nu+1}^{m} P_{\underline{\beta}_m}(A_i) = \alpha,
$$

where the first inequality follows from Lemma 3 (with $\Delta_0 = S_\nu$ and $\Delta_{i-\nu} = \nu X_i/d^*_{\nu,i}$ for $i = \nu+1, \ldots, m$) and the second inequality holds since each term in the summation achieves its maximum at $\underline{\beta}_m$ by Theorem 1. On the other hand, since the thresholds $\{d^*_{\nu,i}\}_{i=\nu+1}^{m-1}$ satisfy $\sum_{i=\nu+1}^{m-1} P_{\underline{\beta}_{m-1}}(A_i) = \alpha$ and each term in this summation satisfies $P_{\underline{\beta}_{m-1}}(A_i) > P_{\underline{\beta}_m}(A_i)$ by Theorem 1, the last term in the summation of (19) is greater than zero, that is, $P_{\underline{\beta}_m}(A_m) > 0$. This guarantees $d^*_{\nu,m}$ to be finite, which implies that rejection region $R^*_{\nu,m}$ is larger than $R^*_{\nu,m-1}$.

Finally, for $m = k$, since the null hypothesis $H_{0,k}$ now contains only one parameter configuration $\underline{\beta}_k$ and $d^*_{\nu,i}$ is available up to $i = k - 1$, one determines $d^*_{\nu,k}$ by solving

$$P_{\underline{\beta}_k}(R^*_{\nu,k}) = \alpha, \tag{21}$$

which implies that $R^*_{\nu,k}$ is level-$\alpha$. Similarly, one can show that $d^*_{\nu,k}$ is finite. Thus, $R^*_{\nu,k}$ is larger than $R^*_{\nu,k-1}$. The determination of $\{d^*_{\nu,m}\}_{m=\nu+1}^k$ is completed.

To conduct the simultaneous tests for $\mathcal{B}$, assert not $H_{0,m}$ (i.e., assert $H_{A,m}$)

$$\text{if } R^*_{\nu,m} \text{ is true.} \tag{22}$$

Notice two facts: (1) $\mathcal{B}$ is closed under the operation of intersection and (2) for each $\nu + 1 \leq m \leq k$, $R^*_{\nu,m} = \bigcap_{i=m}^{k} R^*_{\nu,i}$ is level-$\alpha$. Therefore, the experimentwise error rate is no greater than $\alpha$ by the closed test procedure.

The discussion above is now summarized as the following theorem.

THEOREM 3. *The rejection regions $R^*_{\nu,m}$ given in (17) increase when $m$ gets larger and each defines a level-$\alpha$ test for $H_{0,m}$. If one conducts simultaneous tests for $\mathcal{B}$ using (22), then the experimentwise error rate is controlled at $\alpha$ in the strong sense.*

Let $[1], \ldots, [k]$ be random indices such that $\hat{\beta}^2_{[1]} < \cdots < \hat{\beta}^2_{[k]}$. We now describe the step-up testing procedure based on $R^*_{\nu,m}$ as follows:

*Step* 1: If $R^*_{\nu,\nu+1}$ is true, then conclude that $\beta_{[\nu+1]}, \ldots, \beta_{[k]}$ are the $k - \nu$ active effects ($= H \cap H_{A,\nu+1}$) and stop; otherwise, go to step 2.

*Step* 2: If $R^*_{\nu,\nu+2}$ is true, then conclude that $\beta_{[\nu+2]}, \ldots, \beta_{[k]}$ are the $k - \nu - 1$ active effects ($= H_{0,\nu+1} \cap H_{A,\nu+2}$) and stop; otherwise, go to step 3.

⋮

*Step* $k - \nu$: If $R^*_{\nu,k}$ is true, then conclude that $\beta_{[k]}$ is the only active effect ($= H_{0,k-1} \cap H_{0,k}$) and stop; otherwise, conclude no active effect and stop.

4.2. *The construction of $\{R^*_{m-1,m}\}_{m=\nu+1}^k$: step-up tests with sequential scaling (SUS).* There is another way to conduct the simultaneous tests for $\mathcal{B}$. For each integer $m \in [\nu+1, k]$, we construct a level-$\alpha$ region for $H_{0,m}$, denoted by $R^*_{m-1,m}$ (corresponding to $R_{m-1,m}$ in Section 3), of the form

$$(23) \quad \begin{aligned} R^*_{m-1,m} &= \bigcup_{i=\nu+1}^{m} \{W_{i-1,i} > d^*_{i-1,i}\} = \bigcup_{i=\nu+1}^{m} \{S_\nu < Q_i\} \\ &= \{S_\nu < \max\{Q_i\}_{i=\nu+1}^{m}\}, \end{aligned}$$

where $Q_i = (i-1)X_i/d^*_{i-1,i} - S_{i-1} + S_\nu$ and $S_{i-1} = \sum_{j=1}^{i-1} X_j$.

To determine constants $\{d^*_{m-1,m}\}_{m=\nu+1}^{k}$, we first let the constant $d^*_{\nu,\nu+1}$ equal $d_{\nu,\nu+1}$. Suppose that $d^*_{i-1,i}$ is available up to $i = m-1$ for $m < k$. We then determine $d^*_{m-1,m}$. Comparing (23) with (17), $R^*_{m-1,m}$ and $R^*_{\nu,m}$ have similar forms. Therefore, similarly to (19), we obtain $d^*_{m-1,m}$ by solving

$$(24) \quad \sum_{i=\nu+1}^{m} P_{\underline{\beta}_m}(\max\{S_\nu, \{Q_j\}_{j=\nu}^{i-1}\} < Q_i) = \alpha$$

(with $Q_\nu = 0$). Finally, for $m = k$, since $d^*_{m-1,m}$ is available up to $m = k-1$, $d^*_{k-1,k}$ is solved by $P_{\underline{\beta}_k}(R^*_{k-1,k}) = \alpha$. The determination of $\{R^*_{m-1,m}\}_{m=\nu+1}^{k}$ is thus complete.

Using a discussion similar to that used for $R^*_{\nu,m}$, one can show that $R^*_{m-1,m}$ is a level-$\alpha$ test for $H_{0,m}$ and is increasing in $m$. More specifically, Lemma 3 implies that

$$P_{\underline{\beta}}(R^*_{m-1,m})$$
$$\leq \sum_{i=\nu+1}^{m} P_{\underline{\beta}}(\max\{S_\nu, \{Q_j\}_{j=\nu}^{i-1}\} < Q_i)$$
$$= \sum_{i=\nu+1}^{m} P_{\underline{\beta}}\left(\max\left\{S_{i-1}, \left\{\frac{(j-1)X_j}{d^*_{j-1,j}} + S_{i-1} - S_{j-1}\right\}_{j=\nu}^{i-1}\right\} < \frac{(i-1)X_i}{d^*_{i-1,i}}\right).$$

The last step rewrites each set and makes it clear that, by Theorem 1, each probability above on the right-hand side achieves its maximum at $\underline{\beta}_m$ among $\underline{\beta} \in H_{0,m}$. Therefore, the type I error of $R^*_{m-1,m}$ is bounded by $\alpha$ due to (24). Again, by Theorem 1, each term corresponding to $i < m$ evaluated at $\underline{\beta}_m$ is smaller than that at $\underline{\beta}_{m-1}$, which ensures the existence of a finite solution for $d^*_{m-1,m}$.

To conduct the simultaneous tests for the null hypotheses in $\mathcal{B}$, assert not $H_{0,m}$ (i.e., assert $H_{A,m}$)

(25) \qquad if $R^*_{m-1,m}$ is true.

Therefore, we have a theorem similar to Theorem 3.

THEOREM 4. *The rejection regions $R^*_{m-1,m}$ given in (23) are increasing when $m$ increases and each defines a level-$\alpha$ test for $H_{0,m}$. If one conducts the simultaneous tests for $\mathcal{B}$ using (25), then the experimentwise error rate is strongly controlled at $\alpha$.*

We omit the description of the step-up testing procedure based on $R^*_{m-1,m}$.

REMARK 1. Langsrud and Næs [7] and Venter and Steel [11] also considered these two step-up procedures. For the same test statistics $W_{\nu,m}$ and $W_{m-1,m}$, they proposed to determine critical values $d^\dagger_{\nu,m}$ and $d^\dagger_{m-1,m}$ iteratively by

$$P_{\underline{\beta}_m}\left(\bigcup_{i=\nu+1}^{m}\{W_{\nu,i} > d^\dagger_{\nu,i}\}\right) = \alpha \quad \text{and} \quad P_{\underline{\beta}_m}\left(\bigcup_{i=\nu+1}^{m}\{W_{i-1,i} > d^\dagger_{i-1,i}\}\right) = \alpha.$$
(26)

Intuitively, the solutions $d^\dagger_{\nu,m}$ and $d^\dagger_{m-1,m}$ to the above equations would be smaller than their corresponding cutoff points $d^*_{\nu,m}$ and $d^*_{m-1,m}$ determined by (19) and (24) and would hence result in larger rejection regions. However, it is still not clear that the error rates of their procedures are controlled at $\alpha$ in the strong sense because it is very difficult to establish that for all $\underline{\beta} \in H_{0,m}$,

(27) $\quad P_{\underline{\beta}}(R^*_{\nu,m}) \le P_{\underline{\beta}_m}(R^*_{\nu,m}) \quad \text{or} \quad P_{\underline{\beta}}(R^*_{m-1,m}) \le P_{\underline{\beta}_m}(R^*_{m-1,m}).$

If, for example, we write $R^*_{\nu,m}$ in the form

$$T_L(X_1,\ldots,X_\nu)(=:S_\nu) < T_R(X_{\nu+1},\ldots,X_m)(=:\max\{\nu X_i/d^*_{\nu,i}\}_{i=\nu+1}^m),$$

then $T_R$ involves more than one argument and Theorem 1 cannot be applied. However, our numerical studies show no evidence against (27).

4.3. *An example.* We illustrate the proposed methods using a $2^4$ factorial experiment from [10], pages 246–254, which investigates how temperature, pressure, concentration of formaldehyde and stirring rate influence the filtration rate of a chemical product. The results are presented in Table 1. Column 2 of Table 1 gives the eight effect estimates with largest absolute values and Column 3 the corresponding squared statistics, while $S_7 = \sum_{i=1}^{7} X_i$ equals 15.11 for the seven effect estimates with smallest absolute values. Test statistics $W_{\nu,m}$ and $W_{m-1,m}$ are presented in the next two columns for $\nu = 7$. The SUF procedure identifies four largest active effects, irrespective of the two ways of choosing thresholds ($d^*_{\nu,m}$ or $d^\dagger_{\nu,m}$), while the SUS procedure identifies five largest active effects, also irrespective of the two threshold selections. For the sake of comparison, a step-down procedure from Voss and Wang [14], which uses test statistics $T_{\text{SD},m} = X_m/\min\{0.92S_7, 0.23S_{11}\}$, identifies three largest active effects. Finally, a MATLAB program for the evaluation of cutoff points is available from the authors.

TABLE 1
*The Montgomery* [10] *data, the step-up tests with fixed scaling* (*SUF*) *and the step-up tests with sequential scaling* (*SUS*), *and the related cutoff points*

| Effect | | | Test statistics | | The cutoff points at level $\alpha = 0.05$ | | | |
|---|---|---|---|---|---|---|---|---|
| $m$ | estimate | $X_m$ | $W_{\nu,m}$ | $W_{m-1,m}$ | $d^{\dagger}_{\nu,m}$ | $d^{*}_{\nu,m}$ | $d^{\dagger}_{m-1,m}$ | $d^{*}_{m-1,m}$ |
| 8 | $-2.625$ | 6.89 | 3.2 | 3.2 | 14.9 | 14.9 | 14.9 | 14.9 |
| 9 | 3.125 | 9.77 | 4.5 | 3.6 | 26.5 | 28.0 | 16.4 | 16.7 |
| 10 | 4.125 | 17.02 | 7.9 | 4.8 | 38.4 | 42.0 | 16.0 | 16.3 |
| 11 | 9.875 | 97.52 | 45.2 | <u>20.0</u> | 52.2 | 58.5 | 15.5 | 15.7 |
| 12 | 14.625 | 213.89 | <u>99.1</u> | 16.1 | 67.7 | 77.5 | 15.1 | 15.2 |
| 13 | 16.625 | 276.39 | 128.0 | 9.2 | 85.0 | 99.1 | 14.6 | 14.8 |
| 14 | $-18.125$ | 328.52 | 152.2 | 6.7 | 104.5 | 124.1 | 14.3 | 14.5 |
| 15 | 21.625 | 467.64 | 216.7 | 6.8 | 126.3 | 123.4 | 14.0 | 13.9 |

**5. Simulation study.** A limited simulation study was conducted to compare five testing procedures: step-up tests with sequential scaling (SUS and SUSI using cutoff points determined by (24) and (26), respectively), step-up tests with fixed scaling (SUF and SUFI using cutoff points determined by (19) and (26), respectively) and the Voss and Wang [14] step-down procedure (SD). The testing procedures were evaluated in terms of four measures: (1) the experimentwise error rate (EER), (2) the probability of correctly selecting the number of inactive effects (PCSN), (3) the probability of complete correct selection (PCCS) and (4) the expected fraction of active effects that are declared active (Power). The simulation was carried out for several choices of $k$. Since the results are similar, we only present the choice $k = 15$ on six cases below, following Venter and Steel [11]:

C1: $\underline{\boldsymbol{\beta}} \in H_{0,14}$, $\beta_{15} = s$;  
C2: $\underline{\boldsymbol{\beta}} \in H_{0,12}$; $\beta_{13} = \beta_{14} = \beta_{15} = s$;  
C3: $\underline{\boldsymbol{\beta}} \in H_{0,10}$, $\beta_{11} = \cdots = \beta_{15} = s$;  
C4: $\underline{\boldsymbol{\beta}} \in H_{0,8}$, $\beta_9 = \cdots = \beta_{15} = s$;  
C5: $\underline{\boldsymbol{\beta}} \in H_{0,12}$, $\beta_{12+i} = is, 1 \le i \le 3$;  
C6: $\underline{\boldsymbol{\beta}} \in H_{0,10}$, $\beta_{10+i} = is, 1 \le i \le 5$,

where $s$ takes values from 0 to 8 with a step of 0.02. Each independent sample consists of 15 observations, each from $N(\beta_i, 1)$ for $1 \le i \le 15$.

Simulation results for PCCS are nearly the same as those for PCSN and hence are not reported. Figure 1 presents selected results for the other three evaluation measures for $\alpha = 0.05$ and $\nu = 7$, although findings are similar for other choices of $\alpha$ and $\nu$. Each point was determined based on 100,000 simulations. In summary, all procedures control the EER. Second, there is a very small difference between the two ways of choosing cutoff points, especially between SUS and SUSI. Third, in C1, C2 and C5, the SUS is clearly the best. In C4 and C6, there is a small difference between SUF and SUS. In C3, the SUF performs better at small $s$, but the SUS is better at large $s$. The SD seems to be the worst in most selected cases.

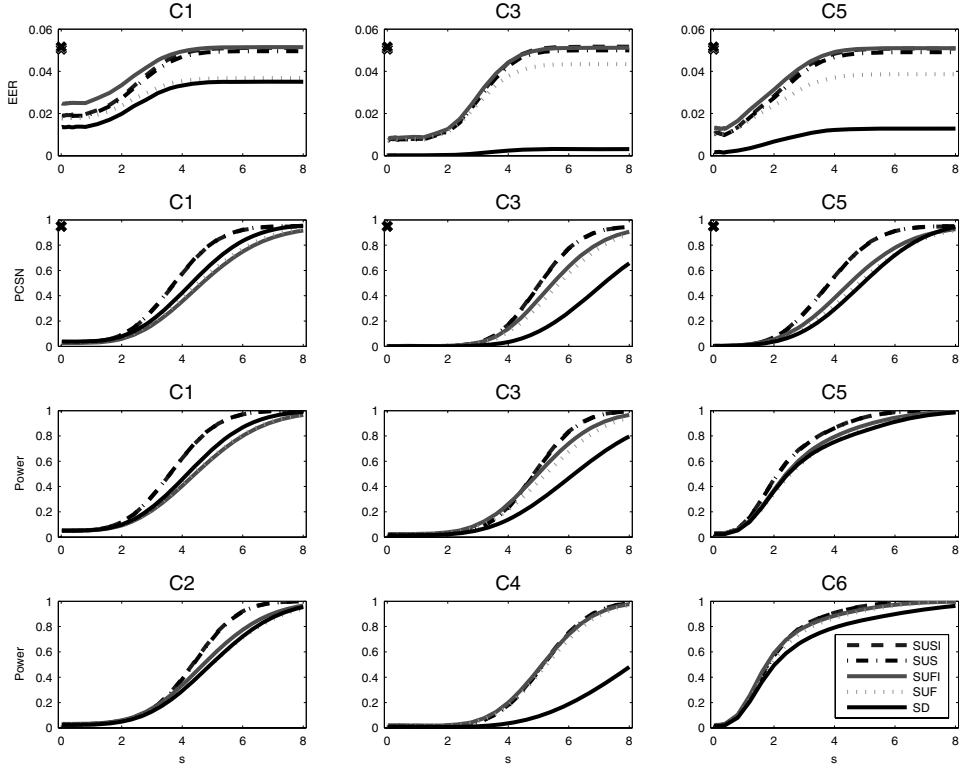

FIG. 1. *Selected simulation results for five test procedures (SUSI, SUS, SUFI, SUF and SD) using $\alpha = 0.05, \nu = 7$, three evaluation measures (EER, PCSN, Power) and six cases of parameter configuration (C1–C6) given in the text.*

**6. Discussion.** We search for active effects in orthogonal saturated designs by conducting simultaneous tests on a sequence of decreasing null hypotheses. A general class of level-$\alpha$ tests is provided for testing at least a certain number of active effects and the least favorable distribution is identified to be the one at $\underline{\boldsymbol{\beta}}_m$. Two sets of simultaneous tests are derived with increasing rejection regions and their experimentwise error rates are controlled at $\alpha$ in the strong sense. Between these two sets of tests, the step-up tests with sequential scaling are recommended because our simulation study indicates that $\{R^*_{m-1,m}\}_{m=\nu+1}^k$ has greater power in most cases. Since the maximal type I errors for $\{R^*_{\nu,m}\}_{m=\nu+1}^k$ and $\{R^*_{m-1,m}\}_{m=\nu+1}^k$ at $m = \nu+1$ and $k$ are equal to $\alpha$, simply enlarging the rejection regions cannot yield valid level-$\alpha$ tests.

We can show that Lemma 1 is also true if $U \sim F_{1,n}$, an $F$-distribution with 1 and $n$ degrees of freedom, and $V \sim F_{1,n}(\lambda)$, a noncentral $F$-distribution with 1 and $n$ degrees of freedom and noncentrality parameter $\lambda$. Conse-

STEP-UP TESTS 15

quently, Theorem 1 remains true if we let $X_1,\ldots,X_k$ be the order statistics of independent random variables with noncentral $F$ distributions $F_{1,n}(\lambda_i)$, $1 \leq i \leq k$. This implies that our step-up simultaneous tests also work for squares of independent $t$-statistics.

**Acknowledgments.** The authors are grateful for the efforts and constructive suggestions of two referees and an Associate Editor.

DIVISION OF BIOSTATISTICS
UNIVERSITY OF FLORIDA
GAINESVILLE, FLORIDA 32610
USA
E-MAIL: samwu@biostat.ufl.edu

DEPARTMENT OF MATHEMATICS AND STATISTICS
WRIGHT STATE UNIVERSITY
DAYTON, OHIO 45435
USA
E-MAIL: wwang@math.wright.edu